# A special constant and series with zeta values and harmonic numbers

**Khristo N. Boyadzhiev**

Department of Mathematics and Statistics,

Ohio Northern University, Ada, OH 45810, USA

k-boyadzhiev@onu.edu

**Abstract.** In this paper we demonstrate the importance of a mathematical constant which is the value of several interesting numerical series involving harmonic numbers, zeta values, and logarithms. We also evaluate in closed form a number of numerical and power series.

## 1. Introduction

The purpose of this paper is to demonstrate the importance of the constant

$$M = \int_0^1 \frac{\psi(t+1) + \gamma}{t} \, dt \approx 1.257746$$

where $\psi(z) = \Gamma'(z)/\Gamma(z)$ is the digamma function and $\gamma = -\psi(1)$ is Euler's constant. Theorem 1 in Section 2 shows that $M$ is the numerical value of many interesting and important series involving logarithms, harmonic numbers $H_n$ and zeta values $\zeta(n)$. The harmonic numbers are defined by

$$H_n = \psi(n+1) + \gamma = 1 + \frac{1}{2} + \ldots + \frac{1}{n}, \ n \geq 1, \quad H_0 = 0$$

with generating function

$$-\frac{\ln(1-t)}{1-t} = \sum_{n=1}^{\infty} H_n t^n, \quad |t|<1.$$

We also use the skew-harmonic numbers given by

$$H_n^- = 1 - \frac{1}{2} + \frac{1}{3} + \ldots + \frac{(-1)^{n-1}}{n}, \quad n \geq 1. \quad H_0^- = 0.$$

Riemann's zeta function is defined by

$$\zeta(s) = \sum_{n=1}^{\infty} \frac{1}{n^s}, \quad \operatorname{Re} s > 1.$$

We consider series with terms $H_n(\zeta(n)-1)$, $H_n(\zeta(n)-\zeta(n+1))$, $H_n^-\big(\zeta(n)-\zeta(n+1)\big)$ and others similar. Related power series with coefficients $H_n \zeta(n+1,a)$ and $H_n^- \zeta(n+1,a)$, where $\zeta(s,a)$ is the Hurwitz zeta function are also studied. Some series are evaluated in closed form in terms of known constants and special functions. Other interesting series will be presented in integral form. The main results of this paper are Theorem 1 in Section 2, Theorem 2 in Section 3, and Theorem 3 in Section 4.

## 2. A mosaic of series

A classical theorem rooted in the works of Christian Goldbach (1690-1764) says that

$$\sum_{n=2}^{\infty} \{\zeta(n)-1\} = 1$$

(see [1], [9], and [12, p. 142]; a short proof is given in the Appendix). The terms $\zeta(n)-1$ of the above series decrease like the powers of $\frac{1}{2}$. The following lemma was proved in [7, p. 51]. It will be needed throughout the paper. (For convenience a proof is given in the Appendix.)

**Lemma 1**. *For every $n \geq 2$*

$$\frac{1}{2^n} < \zeta(n) - 1 < \frac{1}{2^n}\frac{n+1}{n-1}.$$

The lemma shows that the power series with general term $\{\zeta(n)-1\}x^n$ is convergent in the disk $|x|<2$. The following representation is true [12, p. 173, (139)].

(1) $$\sum_{n=2}^{\infty}\{\zeta(n)-1\}x^{n-1} = 1-\gamma-\psi(2-x).$$

Integration of this series yields [12, p. 173, (135)]

(2) $$\sum_{n=2}^{\infty}\frac{\zeta(n)-1}{n}x^n = (1-\gamma)x + \ln\Gamma(2-x)$$

and in particular, when $x=1$, this turns into a classical result published by Euler in 1776 (see [11])

(3) $$\sum_{n=2}^{\infty}\frac{\zeta(n)-1}{n} = 1-\gamma.$$

With $x=-1$ we find

$$\sum_{n=2}^{\infty}(-1)^n\frac{\zeta(n)-1}{n} = \sum_{n=2}^{\infty}(-1)^n\frac{\zeta(n)}{n} - \sum_{n=2}^{\infty}\frac{(-1)^n}{n} = \gamma - 1 + \ln\Gamma(3)$$

so that (as $\ln 2 = \ln\Gamma(3)$)

(4) $$\sum_{n=2}^{\infty}(-1)^n\frac{\zeta(n)}{n} = \gamma$$

which is another representation of $\gamma$ often used by Euler.

Cases when the series can be evaluated in explicit closed form are intertwined with similar cases when this is not possible (unless we accept integral form). Series of the form

$$\sum_{n=2}^{\infty}\frac{\zeta(n)-1}{n+p}, \quad p=1,2$$

have been studied and evaluated in a number of papers - see, for instance, [9], [12, pp. 213-219, (469), (474), 517), (518)] and also [13], [14]. At the same time the series

$$\sum_{n=2}^{\infty}\frac{\zeta(n+1)-1}{n}, \quad \sum_{n=2}^{\infty}(-1)^{n-1}\frac{\zeta(n+1)-1}{n}$$

resists evaluation in closed form in terms of recognized constants. However, they can be evaluated in integral form and we want to point out one simple integral which appears in several interesting cases as evidenced by Theorem 1 below. Namely, consider the constant

$$M = \int_0^1 \frac{\psi(t+1)+\gamma}{t}\,dt.$$

**Theorem 1**. *With $M$ as defined above*,

(a) $$\sum_{n=1}^{\infty}\frac{(-1)^{n-1}\zeta(n+1)}{n} = M$$

(b) $$\sum_{n=1}^{\infty}\frac{1}{n}\ln\left(1+\frac{1}{n}\right) = M$$

(c) $$\sum_{n=1}^{\infty}\frac{\ln(n+1)}{n(n+1)} = M$$

(d) $$\sum_{n=1}^{\infty}H_n\left(\zeta(n+1)-1\right) = M$$

(e) $$\sum_{n=2}^{\infty}H_n\left(\zeta(n)-1\right) = M+1-\gamma$$

(f) $$\sum_{n=1}^{\infty}\frac{1}{n}\left(n-\zeta(2)-\zeta(3)-\ldots-\zeta(n)\right) = M$$

(*the first term in this series is just* 1)

(g) $$\sum_{n=1}^{\infty}H_n^{-}\left(\zeta(n+1)-\zeta(n+2)\right) = M-\ln 2$$

(h) $$\int_0^\infty \frac{\mathrm{Ein}(x)}{e^x - 1} dx = M$$

*where* $\mathrm{Ein}(x)$ *is the modified exponential integral, an entire function*

$$\mathrm{Ein}(x) = \sum_{n=1}^\infty \frac{(-1)^{n-1} x^n}{n! \, n}$$

(i) $$\int_0^1 \frac{(1-u)\ln(1-u)}{u \ln u} du = M$$

(j) $$\sum_{n=1}^\infty \frac{1}{n 2^n} \left\{ \sum_{k=1}^n \binom{n}{k} (-1)^{k-1} \zeta(k+1) \right\} = M.$$

Before proving the theorem we want to mention that the evaluations (a), (b), (c), and (i) are not new (for instance, they are listed on p. 142 in [3]). The series (c) appears in important applications in number theory – see [5] and [8]. The entire paper [4] is dedicated to that series.

*Proof.* Starting from the well-known Taylor expansion for $\psi(1+x) + \gamma$

(5) $$\sum_{n=1}^\infty (-1)^{n-1} \zeta(n+1) x^n = \psi(1+x) + \gamma$$

we divide both sides by $x$ and integrate from $0$ to $1$ to prove (a). Then (a) implies (b)

$$\sum_{n=1}^\infty \frac{(-1)^{n-1} \zeta(n+1)}{n} = \sum_{n=1}^\infty \frac{(-1)^{n-1}}{n} \left\{ \sum_{k=1}^\infty \frac{1}{k^{n+1}} \right\} = \sum_{k=1}^\infty \frac{1}{k} \left\{ \sum_{n=1}^\infty \frac{(-1)^{n-1}}{n} \left( \frac{1}{k} \right)^n \right\} = \sum_{k=1}^\infty \frac{1}{k} \ln\left(1 + \frac{1}{k}\right).$$

Next we prove (b) → (d)

$$\sum_{n=1}^\infty H_n \left( \zeta(n+1) - 1 \right) = \sum_{n=1}^\infty H_n \left\{ \sum_{k=2}^\infty \frac{1}{k^{n+1}} \right\} = \sum_{k=2}^\infty \frac{1}{k} \left\{ \sum_{n=1}^\infty H_n \left( \frac{1}{k} \right)^n \right\}$$

and using the generating functions for the harmonic numbers we continue this way

$$= \sum_{k=2}^{\infty} \frac{1}{k} \left\{ \frac{-1}{1-\frac{1}{k}} \ln\left(1-\frac{1}{k}\right) \right\} = \sum_{k=2}^{\infty} \frac{1}{k-1} \ln\left(\frac{k}{k-1}\right) = \sum_{m=1}^{\infty} \frac{1}{m} \ln\left(\frac{m+1}{m}\right) = M.$$

Now starting from the first sum in the last equation we write

$$M = \sum_{k=2}^{\infty} \frac{1}{k} \left\{ \frac{-1}{1-\frac{1}{k}} \ln\left(1-\frac{1}{k}\right) \right\} = \sum_{k=2}^{\infty} \frac{1}{k-1} \left\{ -\ln\left(1-\frac{1}{k}\right) \right\} = \sum_{k=2}^{\infty} \frac{1}{k-1} \left\{ \sum_{m=1}^{\infty} \frac{1}{mk^m} \right\}$$

$$= \sum_{n=1}^{\infty} \frac{1}{n} \left\{ \sum_{k=2}^{\infty} \frac{1}{k^n (k-1)} \right\} = \sum_{n=1}^{\infty} \frac{1}{n} \left( n - \zeta(2) - \zeta(3) - \ldots - \zeta(n) \right)$$

by using the evaluation (see [7, p. 51])

$$\sum_{k=2}^{\infty} \frac{1}{k^n (k-1)} = n - \zeta(2) - \zeta(3) - \ldots - \zeta(n).$$

Thus (f) is also proved. Equation (e) follows from (d) by writing

$$\sum_{n=2}^{\infty} H_n \left( \zeta(n) - 1 \right) = \sum_{n=2}^{\infty} \left( H_{n-1} + \frac{1}{n} \right) \left( \zeta(n) - 1 \right) = \sum_{m=1}^{\infty} H_m \left( \zeta(m+1) - 1 \right) + \sum_{n=2}^{\infty} \frac{1}{n} \left( \zeta(n) - 1 \right)$$

and then applying (3).

Next we prove (h). From the representation

$$\zeta(s) = \frac{1}{\Gamma(s)} \int_0^{\infty} \frac{x^{s-1}}{e^x - 1} dx, \quad \mathrm{Re}(s) > 1$$

we find with $s = n+1$

$$\sum_{n=1}^{\infty} \frac{(-1)^{n-1} \zeta(n+1)}{n} = \int_0^{\infty} \left\{ \sum_{n=1}^{\infty} \frac{(-1)^{n-1} x^n}{n! n} \right\} \frac{dx}{e^x - 1} = \int_0^{\infty} \frac{\mathrm{Ein}(x)}{e^x - 1} dx.$$

The exponential integral has the representation

$$\operatorname{Ein}(x) = \int_0^x \frac{1-e^{-t}}{t}\, dt$$

and integrating by parts we find

$$\int_0^\infty \frac{\operatorname{Ein}(x)}{e^x - 1}\, dx = \int_0^\infty \frac{\operatorname{Ein}(x)\, d(1-e^{-x})}{1-e^{-x}} = \operatorname{Ein}(x)\ln(1-e^{-x})\Big|_0^\infty - \int_0^\infty \operatorname{Ein}'(x)\ln(1-e^{-x})\, dx$$

$$= -\int_0^\infty \frac{1-e^{-x}}{x}\ln(1-e^{-x})\, dx\ .$$

Here $\operatorname{Ein}(x)\ln(1-e^{-x})\Big|_0^\infty = 0$ since $\operatorname{Ein}(x)$ grows like $\ln x$ when $x \to \infty$ and also $\lim_{x \to 0} x\ln(1-e^{-x}) = 0$. In the last integral we make the substitution $u = e^{-x}$ to find

$$-\int_0^\infty \frac{1-e^{-x}}{x}\ln(1-e^{-x})\, dx = \int_0^1 \frac{(1-u)\ln(1-u)}{u\ln u}\, du$$

and (i) is proved.

We shall prove now the implication (a) → (g) by using Abel's lemma for transformation of series (see, for example, [10, Exercise 10, p.78]):

**Lemma 2.** *Let $\{a_n\}$ and $\{b_n\}$, $n \geq p$, be two sequences of complex numbers and let $A_n = a_p + a_{p+1} + \ldots + a_n$. Then for every $n > p$ we have*

$$\sum_{k=p}^n a_k b_k = b_n A_n + \sum_{k=p}^{n-1} A_k (b_k - b_{k+1})\ .$$

We take here $p = 1$, $a_k = \dfrac{(-1)^{k-1}}{k}$, and $b_k = \zeta(k+1) - 1$. Then $A_k = H_k^-$ and we find

$$\sum_{k=1}^n \frac{(-1)^{k-1}}{k}\{\zeta(k+1) - 1\} = H_n^-\{\zeta(n+1) - 1\} + \sum_{k=1}^{n-1} H_k^-\{\zeta(k+1) - \zeta(k+2)\}.$$

Setting $n \to \infty$ and using the estimate from lemma1 and also the fact that $|H_n^-| \le H_n$ and $H_n \square \ln n$ at infinity we come to the equation

$$\sum_{k=1}^{\infty} \frac{(-1)^{k-1}}{k}\{\zeta(k+1)-1\} = \sum_{k=1}^{\infty} H_k^- \{\zeta(k+1)-\zeta(k+2)\}.$$

According to (a) we have

$$\sum_{k=1}^{\infty} \frac{(-1)^{k-1}}{k}\{\zeta(k+1)-1\} = M - \ln 2$$

and (g) follows.

The implication (b) $\to$ (c) also follows from Abel's lemma. We take $p=1$ and

$$a_k = \ln\left(1+\frac{1}{k}\right) = \ln(k+1)-\ln(k), \quad b_k = \frac{1}{k}.$$

Here $A_n = \ln(n+1)$ and lemma 2 provides the equation

$$\sum_{k=1}^{n} \frac{1}{k}\ln\left(1+\frac{1}{k}\right) = \frac{\ln(n+1)}{n} + \sum_{k=1}^{n-1} \frac{\ln(k+1)}{k(k+1)}$$

which clearly shows the relation between (b) and (c).

Last we prove (j). The proof is based on Euler's series transformation [2]. Given a power series $f(x) = a_0 + a_1 x + \ldots$ we have for sufficiently small $|t|$

$$\frac{1}{1-t} f\left(\frac{t}{1-t}\right) = \sum_{n=0}^{\infty} t^n \left\{ \sum_{k=0}^{n} \binom{n}{k} a_k \right\}.$$

We take $f(x) = \psi(x+1) + \gamma$ with the expansion (5), where $a_0 = f(0) = 0$. Using the substitution $x = \frac{t}{1-t}$ we compute

$$\int_0^1 \frac{\psi(x+1)+\gamma}{x} dx = \int_0^{1/2} \frac{1}{1-t} f\left(\frac{t}{1-t}\right) \frac{dt}{t}$$

$$= \int_0^{1/2} \sum_{n=1}^\infty t^n \left\{ \sum_{k=1}^n \binom{n}{k} (-1)^{k-1} \zeta(k+1) \right\} \frac{dt}{t} = \sum_{n=1}^\infty \frac{1}{n 2^n} \sum_{k=1}^n \binom{n}{k} (-1)^{k-1} \zeta(k+1)$$

after integrating term by term.

It is easy to see that for $n \geq 1$

$$\sum_{k=1}^n \binom{n}{k} (-1)^{k-1} \zeta(k+1) = \int_0^\infty \frac{1 - L_n(x)}{e^x - 1} dx$$

where $L_n(x)$ are the Laguerre polynomials.

The proof of the theorem is completed.

**Remark 1**. Applying Lemma 2 for $p = 1$, $a_k = \frac{1}{k}$, and $b_k = \zeta(k+1) - 1$ we find

$$\sum_{k=1}^\infty \frac{\zeta(k+1) - 1}{k} = \sum_{k=1}^\infty H_n \{ \zeta(k+1) - \zeta(k+2) \}.$$

At the same time from (1)

$$\sum_{k=1}^\infty \frac{\zeta(k+1) - 1}{k} = \int_0^1 \frac{1 - \gamma - \psi(2-t)}{t} dt \ .$$

**Remark 2**. This is a note on the series in (f). Let

$$S_n = \frac{1}{2^n} + \frac{1}{3^n 2} + \frac{1}{4^n 3} + \ldots = \sum_{k=2}^\infty \frac{1}{k^n (k-1)} = n - \zeta(2) - \zeta(3) - \ldots - \zeta(n).$$

Obviously,

$$\frac{1}{2^n} < S_n < \frac{1}{2^n} + \frac{1}{3^n} + \frac{1}{4^n} + \ldots = \zeta(n) - 1 < \frac{1}{2^n} \frac{n+1}{n-1}.$$

That is, we have the same estimate as in Lemma 1

$$\frac{1}{2^n} < S_n < \frac{1}{2^n} \frac{n+1}{n-1} \ .$$

It is interesting to compare equation (g) from the above theorem with the following result.

**Proposition 1.** *We have*

(6) $$\sum_{n=2}^{\infty} H_n^- \left(\zeta(n)-\zeta(n+1)\right) = \frac{\pi^2}{6} - \gamma - \ln 2.$$

The proof comes from Abel's lemma again, where we take $p=2$, $a_k = \dfrac{(-1)^{k-1}}{k}$ and $b_k = \zeta(k)-1$. Then $A_k = H_k^- - 1$ and the lemma yields for every $n \geq 3$

$$\sum_{k=2}^{n} \frac{(-1)^{k-1}}{k}\{\zeta(k)-1\} = (H_n^- - 1)\{\zeta(n)-1\} + \sum_{k=2}^{n-1}(H_k^- -1)\{\zeta(k)-\zeta(k+1)\}.$$

From here with $n \to \infty$ we find

$$\sum_{k=2}^{\infty} \frac{(-1)^{k-1}}{k}\zeta(k) - \ln 2 + 1 = \sum_{k=2}^{\infty} H_k^-\{\zeta(k)-\zeta(k+1)\} - \sum_{k=2}^{\infty}\{\zeta(k)-\zeta(k+1)\}.$$

The first sum equals $-\gamma$ and for the last sum we write

$$\sum_{k=2}^{\infty}\{\zeta(k)-\zeta(k+1)\} = \sum_{k=2}^{\infty}\{(\zeta(k)-1) - (\zeta(k+1)-1)\}$$

which is a telescoping series equal to its first term $\zeta(2)-1$. Thus

$$-\gamma - \ln 2 + 1 = \sum_{k=2}^{\infty} H_k^-\{\zeta(k)-\zeta(k+1)\} - \zeta(2) + 1.$$

and the proof is finished.

It was noted by the referee that another proof follows from the fact that the series in (6) has a telescoping property. This property brings to a convenient representation of its partial sum

$$\sum_{n=2}^{m} H_n^- \left(\zeta(n)-\zeta(n+1)\right) = \frac{\zeta(2)}{2} - H_{m+1}^-(\zeta(m+1)-1) - H_{m+1}^- + \sum_{n=2}^{m} \frac{(-1)^n \zeta(n+1)}{n+1}.$$

Setting here $m \to \infty$ and using equation (4) we come the desired evaluation.

It is also interesting to compare parts (a) and (c) of Theorem 1 to the result in the following proposition.

**Proposition 2.** *For any $p > 1$*

$$\sum_{n=1}^{\infty} \frac{\log(n+1)}{n^p} = -\zeta'(p) + \sum_{k=1}^{\infty} \frac{(-1)^{k-1} \zeta(p+k)}{k}.$$

*Proof.*

$$\sum_{n=1}^{\infty} \frac{\log(1+n)}{n^p} = \sum_{n=1}^{\infty} \frac{1}{n^p} \log\left[n\left(1+\frac{1}{n}\right)\right] = \sum_{n=1}^{\infty} \frac{\log n}{n^p} + \sum_{n=1}^{\infty} \frac{1}{n^p} \log\left(1+\frac{1}{n}\right)$$

$$= -\zeta'(p) + \sum_{n=1}^{\infty} \frac{1}{n^p} \left\{ \sum_{k=1}^{\infty} \frac{(-1)^{k-1}}{k} \left(\frac{1}{n}\right)^k \right\} = -\zeta'(p) + \sum_{k=1}^{\infty} \frac{(-1)^{k-1}}{k} \left\{ \sum_{n=1}^{\infty} \frac{1}{n^{k+p}} \right\}$$

$$= -\zeta'(p) + \sum_{k=1}^{\infty} \frac{(-1)^{k-1} \zeta(k+p)}{k}.$$

The change of order of summation is justified by the absolute convergence of the series.

The next proposition parallels some of the results in Theorem 1.

**Proposition 3**. *Let*

$$M_1 = \int_1^2 \frac{\psi(1+t)+\gamma}{t} dt \approx 0.86062.$$

*Then*

(k) $$\sum_{n=1}^{\infty} \frac{1}{n} \ln\left(1+\frac{1}{n+1}\right) = M_1$$

(l) $$\sum_{n=1}^{\infty} H_n^- \left(\zeta(n+1)-1\right) = M_1$$

(m) $$\sum_{n=1}^{\infty} \frac{(-1)^{n-1}}{n} \left(n - \zeta(2) - \zeta(3) - \ldots - \zeta(n)\right) = M_1$$

(*the first term in the sum is* 1).

Here (k) follows immediately from the well-known representation

$$\frac{\psi(1+t)+\gamma}{t} = \sum_{n=1}^{\infty} \frac{1}{n(n+t)}$$

by integration from 1 to 2 (integration between 0 and 1 implies (b) in Theorem 1). Independently (l) follows from Theorem 3 in Section 4 (see the remark at the end of that section). Further details are left to the reader.

## 3. Variations on a problem of Ovidiu Furdui

The identity (d) of Theorem 1 can be written in the form

$$\sum_{n=2}^{\infty} H_n \left( \zeta(n+1) - 1 \right) = M - \zeta(2) + 1$$

by starting the summation from $n = 2$. Subtracting this from (e) we find

(7) $$\sum_{n=2}^{\infty} H_n \left[ \zeta(n) - \zeta(n+1) \right] = \frac{\pi^2}{6} - \gamma.$$

This result was displayed in Problem W10 from [6]. We shall extend (7) to power series involving the Hurwitz zeta function

$$\zeta(s,a) = \sum_{k=0}^{\infty} \frac{1}{(k+a)^s}, \quad \zeta(s,1) = \zeta(s),$$

where $\operatorname{Re}(s) > 1$ and $a > 0$.

**Theorem 2.** *For $a > 0$ and $|x| < a$ we have*

(8) $$\sum_{n=2}^{\infty} H_n \left[ \zeta(n,a) - x\zeta(n+1,a) \right] x^n = \zeta(2,a)x^2 + \psi(a)x + \log \Gamma(a-x) - \log \Gamma(a).$$

*Proof.* The left hand side can be written this way

$$\sum_{n=2}^{\infty} H_n \zeta(n,a) x^n = \sum_{n=2}^{\infty} \left( H_{n-1} + \frac{1}{n} \right) \zeta(n,a) x^n = \sum_{n=2}^{\infty} H_{n-1} \zeta(n,a) x^n + \sum_{n=2}^{\infty} \frac{1}{n} \zeta(n,a) x^n$$

$$= \sum_{n=1}^{\infty} H_n \zeta(n+1,a) x^{n+1} + \log \Gamma(a-x) - \log \Gamma(a) + \psi(a) x$$

by using the well-known series from [3, p.78], or [12, p. 159]

(9) $$\sum_{n=2}^{\infty} \frac{1}{n} \zeta(n,a) x^n = \log \Gamma(a-x) - \log \Gamma(a) + \psi(a) x$$

for $|x| < a$. Then we separate the first term in the series

$$\sum_{n=1}^{\infty} H_n \zeta(n+1,a) x^{n+1} = \sum_{n=2}^{\infty} H_n \zeta(n+1,a) x^{n+1} + \zeta(2,a) x^2$$

and from here

$$\sum_{n=2}^{\infty} H_n \zeta(n,a) x^n - \sum_{n=2}^{\infty} H_n \zeta(n+1,a) x^{n+1} = \zeta(2,a) x^2 + \psi(a) x + \log \Gamma(a-x) - \log \Gamma(a).$$

Thus (8) is proved.

The convergence in these series is supported by the estimate given in the Appendix

$$\frac{1}{(a+1)^n} + \frac{1}{a^n} < \zeta(n,a) \leq \frac{1}{a^n} + \frac{1}{(a+1)^n} \left( \frac{n+a}{n-1} \right)$$

($n \geq 2$) and by the slow growth of the harmonic numbers at infinity, $H_n \sim \ln n$.

**Corollary 1**. *For $|x| < 1$*

(10) $$\sum_{n=2}^{\infty} H_n [\zeta(n) - x \zeta(n+1)] x^n = \zeta(2) x^2 - \gamma x + \log \Gamma(1-x),$$

and for $|x| < 2$

(11) $$\sum_{n=2}^{\infty} H_n [\zeta(n) - x \zeta(n+1) + x - 1] x^n = \zeta(2) x^2 - x^2 + (1-\gamma) x + \log \Gamma(2-x).$$

*Proof.* Setting $a=1$ in (8) we have $\psi(1)=-\gamma$ and (10) follows. With $a=2$ in (8) we have

$\zeta(n,2)=\zeta(n)-1$, $\zeta(n+1,2)=\zeta(n+1)-1$, and $\psi(2)=1-\gamma$. Thus (8) turns into (11).

With $x=1$ (11) becomes (7). Setting $x=-1$ in (11) we find also

(12) $$\sum_{n=2}^{\infty}(-1)^n H_n\left[\zeta(n)+\zeta(n+1)-2\right]=\zeta(2)+\gamma-2+\log 2.$$

The series in (8) can be regularized, modified to extend the interval of convergence.

**Corollary 2.** *For every $a>0$ and every $|x|<a+1$*

(13) $$\sum_{n=2}^{\infty}H_n\left[\zeta(n,a)-x\zeta(n+1,a)-\frac{a-x}{a^{n+1}}\right]x^n$$

$$=\zeta(2,a)x^2+\psi(a)x+\frac{(a-x)x}{a^2}+\log\Gamma(a+1-x)-\log\Gamma(a+1).$$

*Proof.* For $|x|<a$ we write using the generating function for the harmonic numbers

$$\log(a-x)=(a-x)\frac{\log(a-x)}{a-x}=(a-x)\left\{\left(\frac{1}{a-x}\right)\left(\log a+\log\left(1-\frac{x}{a}\right)\right)\right\}$$

$$=(a-x)\left\{\frac{\log a}{a-x}+\frac{1}{a}\left(1-\frac{x}{a}\right)^{-1}\log\left(1-\frac{x}{a}\right)\right\}$$

$$=\log a-(a-x)\sum_{n=1}^{\infty}\frac{H_n x^n}{a^{n+1}}=\log a-(a-x)\frac{x}{a^2}-\sum_{n=2}^{\infty}H_n\frac{(a-x)}{a^{n+1}}x^n.$$

This result can be put in the form

$$-\sum_{n=2}^{\infty}H_n\frac{(a-x)}{a^{n+1}}x^n=\log(a-x)-\log a+(a-x)\frac{x}{a^2}.$$

Now we add this equation to equation (8). On the left hand side we unite the two sums into one.
On the right hand side we write $\log(a-x)+\log\Gamma(a-x)=\log\Gamma(a+1-x)$ and $\log a + \log \Gamma(a) = \log \Gamma(a+1)$. This way from (8) we obtain (13). The right side in (13) is obviously defined for $|x| < a+1$. To see that the series on the left side converges in this interval we write it in the form

$$\sum_{n=2}^{\infty} H_n \left[ \zeta(n,a) - x\zeta(n+1,a) - \frac{a-x}{a^{n+1}} \right] x^n$$

$$= \sum_{n=2}^{\infty} H_n \left[ \zeta(n,a) - \frac{1}{a^n} \right] x^n - \sum_{n=2}^{\infty} H_n \left[ \zeta(n+1,a) - \frac{1}{a^{n+1}} \right] x^{n+1}.$$

Here both series converge for $|x| < a+1$ in view of the estimate in the Appendix

$$\frac{1}{(a+1)^m} < \zeta(m,a) - \frac{1}{a^m} \leq \frac{1}{(a+1)^m} \left( \frac{m+a}{m-1} \right)$$

true for every integer $m \geq 2$. The proof is completed.

### 4. Two interesting generating functions

The series in (10) is a power series, the difference of the two powers series with terms

$$H_n \zeta(n) x^n \quad \text{and} \quad H_n \zeta(n+1) x^{n+1}.$$

Although the difference of these two series can be evaluated in closed form, evaluating each one separately brings to difficult integrals. We shall demonstrate this in the following theorem. For comparison we also include the series with skew-harmonic numbers $H_n^-$.

**Theorem 3**. *With $a > 0$ and $|x| < a$*

(14) $\quad \displaystyle\sum_{n=1}^{\infty} H_n \zeta(n+1,a) x^n = \sum_{k=1}^{\infty} \frac{(-1)^{k-1} \zeta(k+1, a-x)}{k} x^k = \int_0^x \frac{\psi(a-x+t) - \psi(a-x)}{t} dt$

(15) $\quad \displaystyle\sum_{n=1}^{\infty} H_n^- \zeta(n+1,a) x^n = \sum_{k=1}^{\infty} \frac{(-1)^{k-1}(2^k-1)\zeta(k+1, a-x)}{k} x^k = \int_0^x \frac{\psi(a-x+2t) - \psi(a-x+t)}{t} dt$.

The proof is based on the lemma:

**Lemma 3.** *Let $f(t) = a_1 t + a_2 t^2 + \ldots + a_n t^n + \ldots$ be a power series. Then the following representations hold*

$$(16) \qquad \sum_{n=1}^{\infty} H_n a_n x^n = \sum_{k=1}^{\infty} \frac{(-1)^{k-1}}{k!k} x^k f^{(k)}(x)$$

$$(17) \qquad \sum_{n=1}^{\infty} H_n^- a_n x^n = \sum_{k=1}^{\infty} \frac{(-1)^{k-1}(2^k - 1)}{k!k} x^k f^{(k)}(x).$$

*Proof.* The harmonic numbers have the well-known integral representation

$$H_n = \int_0^1 \frac{t^n - 1}{t - 1} dt.$$

Multiplying both sides here by $a_n x^n$ and summing for $n = 1, 2, \ldots$ we find for $x$ small enough

$$\sum_{n=1}^{\infty} H_n a_n x^n = \int_0^1 \frac{f(xt) - f(x)}{t - 1} dt$$

(cf. [10, Exercise 20, p.79]). Now Taylor's formula applied to the function $g(t) = f(xt)$ and centered at $t = 1$ gives

$$f(xt) = f(x) + \sum_{k=1}^{\infty} \frac{x^k f^{(k)}(x)}{k!} (t - 1)^k.$$

Substituting this into the integral and integrating term by term we arrive at (16).

The proof of (17) is done the same way by using the representations

$$\sum_{n=1}^{\infty} H_n^- a_n x^n = \int_0^1 \frac{f(x) - f(-xt)}{t + 1} dt$$

and

$$f(-xt) = f(x) + \sum_{k=1}^{\infty} \frac{(-1)^k x^k f^{(k)}(x)}{k!} (t + 1)^k.$$

*Proof of the theorem.* We apply the lemma to the function ([12, p.159])

(18) $$f(t) = \sum_{n=1}^{\infty} \zeta(n+1, a)t^n = -\psi(a-t) + \psi(a), \quad |t| < |a|.$$

Here

$$f^{(k)}(x) = (-1)^{k+1}\psi^{(k)}(a-x) = k!\zeta(k+1, a-x)$$

and from the lemma

(19) $$\sum_{n=1}^{\infty} H_n \zeta(n+1, a)x^n = \sum_{k=1}^{\infty} \frac{(-1)^{k-1}\zeta(k+1, a-x)}{k} x^k.$$

Now we write (18) in the form

$$\sum_{k=1}^{\infty} (-1)^{k-1}\zeta(k+1, a)t^k = \psi(a+t) - \psi(a), \quad |t| < |a|.$$

Dividing by $t$ and integrating between $0$ and $x$ we find

$$\sum_{k=1}^{\infty} \frac{(-1)^{k-1}\zeta(k+1, a)}{k} x^k = \int_0^x \frac{\psi(a+t) - \psi(a)}{t} dt.$$

Replacing here $a$ by $a - x$ gives

(20) $$\sum_{k=1}^{\infty} \frac{(-1)^{k-1}\zeta(k+1, a-x)}{k} x^k = \int_0^x \frac{\psi(a-x+t) - \psi(a-x)}{t} dt.$$

Now (19) and (20) bring to (14). The proof of (15) is left to the reader.

**Remark 3.** By setting $a = 2$ and $x = 1$ in (14) we obtain

$$\sum_{n=1}^{\infty} H_n \left[\zeta(n+1) - 1\right] = \sum_{k=1}^{\infty} \frac{(-1)^{k-1}\zeta(k+1)}{k} = \int_0^1 \frac{\psi(1+t) + \gamma}{t} dt$$

that is, a new proof of equations (a) and (d) in Theorem 1, as $\zeta(n+1, 2) = \zeta(n+1) - 1$.

With $a = 2, x = 1$ in (15) we find

(21) $$\sum_{n=1}^{\infty} H_n^{-} \left(\zeta(n+1) - 1\right) = \int_0^1 \frac{\psi(1+2t) - \psi(1+t)}{t} dt = \int_1^2 \frac{\psi(1+t) + \gamma}{t} dt$$

which is (l) in Proposition 3. To show the equality of the above two integrals we write

$$\int_0^1 \frac{\psi(1+2t)-\psi(1+t)}{t}dt = \int_0^1 \frac{\psi(1+2t)+\gamma}{t}dt - \int_0^1 \frac{\psi(1+t)+\gamma}{t}dt$$

and then in the first integral on the right hand side we make the substitution $u = 2t$.

**Appendix**

*Proof of Goldbach's theorem.* The theorem has a simple two-line proof

$$\sum_{n=2}^{\infty}\{\zeta(n)-1\} = \sum_{n=2}^{\infty}\left\{\sum_{k=2}^{\infty}\frac{1}{k^n}\right\} = \sum_{k=2}^{\infty}\left\{\sum_{n=2}^{\infty}\left(\frac{1}{k}\right)^n\right\} = \sum_{k=2}^{\infty}\left\{\frac{1}{k^2}\cdot\frac{1}{1-\frac{1}{k}}\right\}$$

$$=\sum_{k=2}^{\infty}\left\{\frac{1}{k(k-1)}\right\} = \sum_{k=2}^{\infty}\left\{\frac{1}{k-1}-\frac{1}{k}\right\} = 1$$

where the last sum is a well-known telescoping series.

*Proof of Lemma 1* (extended version). Let $a > 0$, $s > 1$. Then

$$\frac{1}{(a+1)^s} < \zeta(s,a) - \frac{1}{a^s} \leq \frac{1}{(a+1)^s}\left(\frac{s+a}{s-1}\right).$$

Furdui's proof from [7, p.51] is modified here for the Hurwitz zeta function. The left hand side inequality is obvious and for the right hand side inequality we write using the reminder estimate from the integral test for series

$$\sum_{k=2}^{\infty}\frac{1}{(k+a)^s} \leq \int_1^{\infty}\frac{1}{(t+a)^s}dt = \frac{1}{(s-1)(a+1)^{s-1}}.$$

Therefore,

$$\zeta(s,a) - \frac{1}{a^s} = \frac{1}{(a+1)^s} + \sum_{k=2}^{\infty}\frac{1}{(k+a)^s} \leq \frac{1}{(a+1)^s} + \frac{1}{(s-1)(a+1)^{s-1}} = \frac{1}{(a+1)^s}\left(\frac{s+a}{s-1}\right).$$

The proof is completed.